\documentclass[11pt]{amsart}
\usepackage{latexsym}
\usepackage{amscd}
\usepackage[all]{xy}
\usepackage[mathscr]{euscript}
\usepackage{amssymb}

\normalsize

\RequirePackage{xspace}

\newcommand{\thought}[1]{}
\renewcommand{\thought}[1]{ \textbf{[#1]}}

\usepackage{enumerate}        
\newenvironment{roenumerate}{\begin{enumerate}[\upshape (i)]}{\end{enumerate}}

\newcommand\nc {\newcommand}
\newcommand\rnc{\renewcommand}
\nc\script{\mathscr}

\newcount\blopone
\newcount\xone
\newcount\xtwo
\newcount\ytwo

\newtheorem{theorem}{Theorem}[section]
\newtheorem{prop}[theorem]{Proposition}
\newtheorem{equivform}[theorem]{Equivalent Formulation}

\newtheorem{refinement}[theorem]{Refinement}
\newtheorem{summary}[theorem]{Summary}
\newtheorem{importnota}[theorem]{Important Notation}
\newtheorem{prblm}[theorem]{Problem}
\newtheorem{notation}[theorem]{Notation}
\newtheorem{defin}[theorem]{Definition}
\newtheorem{caution}[theorem]{Caution}
\newtheorem{remark}[theorem]{Remark}
\newtheorem{reminder}[theorem]{Reminder}
\newtheorem{lemma}[theorem]{Lemma}
\newtheorem{construction}[theorem]{Construction}
\newtheorem{corollary}[theorem]{Corollary}
\newtheorem{example}[theorem]{Example}
\newtheorem{conclusion}[theorem]{Conclusion}
\newtheorem{triviality}[theorem]{Triviality}
\newtheorem{proto}[theorem]{Prototype Quasifibration}
\newtheorem{cauex}[theorem]{Cautionary Example}
\newtheorem{hypo}[theorem]{Hypothesis}
\newtheorem{subth}{ }[theorem]
\newtheorem{case}{Case}[theorem]
\newtheorem{ssubth}{ }[subth]

\nc\tri[1]{\begin{triviality}
\label{#1}}
\nc\eqf[1]{\begin{equivform}
\label{#1}\begin{em}}
\nc\cas[1]{\begin{case}
\label{#1}
\begin{em}}
\nc\rfn[1]{\begin{refinement}
\label{#1}}
\nc\prt[1]{\begin{proto}
\label{#1}}
\nc\lem[1]{\begin{lemma}
\label{#1}}
\nc\pro[1]{\begin{prop}
\label{#1}}
\nc\thm[1]{\begin{theorem}
\label{#1}}
\nc\cor[1]{\begin{corollary}
\label{#1}}
\nc\dfn[1]{\begin{defin}
\begin{em}
\label{#1}}
\nc\sthm[1]{\begin{subth}
\label{#1}}
\nc\exm[1]{\begin{example}
\label{#1}
\begin{em}}
\nc\plm[1]{\begin{prblm}
\label{#1}
\begin{em}}
\nc\rmk[1]{\begin{remark}
\label{#1}
\begin{em}}
\nc\rmd[1]{\begin{reminder}
\label{#1}
\begin{em}}
\nc\ntn[1]{\begin{notation}
\label{#1}
\begin{em}}
\nc\smr[1]{\begin{summary}
\label{#1}
\begin{em}}
\nc\cau[1]{\begin{caution}
\label{#1}
\begin{em}}
\nc\hyp[1]{\begin{hypo}
\label{#1}\begin{em}}
\nc\imn[1]{\begin{importnota}
\label{#1}
\begin{em}}
\nc\cax[1]{\begin{cauex}
\label{#1}
\begin{em}}
\nc\con[1]{\begin{construction}
\label{#1}
\begin{em}}
\nc\ssthm[1]{\begin{ssubth}
\label{#1}
\begin{em}}
\nc\cnc[1]{\begin{conclusion}
\label{#1}
\begin{em}}

\nc\elem{\end{lemma}}
\nc\eeqf{\blimy \end{em}
\end{equivform}}
\nc\erfn{\end{refinement}}
\nc\eprt{\end{proto}}
\nc\ethm{\end{theorem}}
\nc\ecor{\end{corollary}}
\nc\edfn{\blimy \end{em}
\end{defin}}
\nc\esthm{\end{subth}}
\nc\epro{\end{prop}}
\nc\etri{\end{triviality}}
\nc\eexm{\blimy \end{em}
\end{example}}
\nc\ermk{\blimy \end{em}
\end{remark}}
\nc\ermd{\blimy \end{em}
\end{reminder}}
\nc\eplm{\end{em}
\end{prblm}}
\nc\ecas{\end{em}
\end{case}}
\nc\ecau{\end{em}
\end{caution}}
\nc\ecax{\end{em}
\end{cauex}}
\nc\eimn{\end{em}
\end{importnota}}
\nc\entn{\blimy \end{em}
\end{notation}}
\nc\econ{\end{em}
\end{construction}}
\nc\esmr{\end{em}
\end{summary}}
\nc\ehyp{\blimy \end{em}
\end{hypo}}
\nc\ecnc{\end{em}
\end{conclusion}}
\nc\essthm{\end{em}
\end{ssubth}}

\nc\blimy{{\begin{flushright} $\Box$\par
\end{flushright}}\vskip 2mm }

\newtheorem{hypothesis}[theorem]{Hypothesis}
\nc\sst{\scriptstyle}
\newcommand{\comment}[1]{}
\newcommand{\ri}{\longrightarrow}

\nc\bR{{\mathbf R}}
\nc\bS{{\mathbf S}}
\nc\bT{{\mathbf T}}
\nc\bU{{\mathbf U}}
\nc\z{\zeta}
\nc\bc{{\mathbb{BC}}}
\nc\ct{{\script T}}
\nc\cs{{\script S}}
\nc\car{{\script R}}
\nc\ca{{\script A}}
\nc\cb{{\script B}}
\nc\cc{{\script C}}
\nc\cd{{\script D}}
\nc\ck{{\script K}}
\nc\ce{{\script E}}
\nc\ci{{\script I}}
\nc\co{{\script O}}
\nc\cu{{\script U}}
\nc\bZ{{\mathbb Z}}
\nc\bd{\begin{description}}
\nc\ed{\end{description}}
\nc\ctob{{\script C}at\big(\ci^{op},\ca\big)}
\nc\clim{{\ds\mathop{\rm lim}_{\ds\longleftarrow}}}
\nc\climi{\clim^{\!i}\,}
\nc\climn{\clim^{\!n}\,}
\nc\colim{{\ds\mathop{\rm colim}_{\ds\la}}}
\nc\oa{\overline{\ca}}
\nc\s{\sigma}
\nc\ta{\tau}
\nc\os{\overline\sigma}
\nc\ot{\overline\tau}
\nc\T{\Sigma}
\nc\de[1]{{\mathop{\rm deg(#1)}}}
\nc\Ad[1]{\mathop{\rm Ad}(#1)}
\nc\ad[1]{\mathop{\rm ad}(#1)}
\nc\Tor{\text{\rm Tor}}
\nc\be{\begin{roenumerate}}
\nc\ee{\end{roenumerate}}

\def\der #1 {D\left(#1\right)}
\nc\prf{\begin{proof}}
\nc\eprf{\end{proof}}
\nc\ds{\displaystyle}

\nc\ab{{\script A}b}

\nc\csab{{\script C}at\big(\cs^{op},\ab\big)}
\nc\ctab{{\script C}at\Big({\{\ct^\alpha\}}^{op},\ab\Big)}
\nc\csex{{\script E}x\big(\cs^{op},\ab\big)}
\nc\ctex{{\script E}x\Big({\{\ct^\alpha\}}^{op},\ab\Big)}
\nc\sub{\qquad\subset\qquad}
\nc\ctr[1]{{\left.\ct\left(-,#1\right)\right|}_{\cs}}
\nc\ctrf[2]{{\left.\ct\left(#1,#2\right)\right|}_{\cs}}
\nc\Ctr[1]{{\left.\ct\left(-,#1\right)\right|}_{\ct^\alpha}}
\nc\Ctrf[2]{{\left.\ct\left(#1,#2\right)\right|}_{\ct^\alpha}}

\nc\la{\longrightarrow}
\nc\oti{{^L\otimes_A^{}}}
\nc\rs{\s^{-1}A}
\nc\rrs{{\{\s^{-1}A\}}\op}
\nc\br{{\{\s^{-1}A\}}}
\nc\sm{\s^{-1}}

\nc\nin{\noindent}
\nc\cad[1]{\text{card}(#1)}
\nc\eq{\quad=\quad}
\nc\BA{\begin{array}{c}}
\nc\EA{\end{array}}
\nc\kth{{\it K}--theory}

\nc\barr{
\[
\begin{array}{cccccccccccccccc}
}
\nc\earr{
\end{array}
\]
}
\nc\as[1]{{\langle S\rangle}^{#1}}
\nc\sh{\hbox{\it shift}}

\nc\yy[1]{{\left.\ct\left(-,#1\right)\right|}_{\ct^c}}
\nc\vrep[2]{{\left.\ct\left(#1,#2\right)\right|}_{\ct^\alpha}}
\nc\da{\downarrow}
\nc\Hom{{\mathop{\rm Hom}}}
\nc\End{{\mathop{\rm End}}}
\nc\Ext{{\mathop{\rm Ext}}}
\nc\mr\Modtc
\nc\PExt{{\mathop{\rm PExt}}}
\nc\bA{{\mathbf A}}
\nc\bB{{\mathbf B}}
\nc\bC{{\mathbf C}}
\nc\bD{{\mathbf D}}

\nc\ra\ri
\nc\y[1]{\mathbf{y}#1}
\nc\x[1]{\mathbf{z}#1}
\nc\Mod[1]{\ensuremath{\mathop{\textup{Mod-}#1}}\xspace}
\nc\Md {\ensuremath{\mathop{\textup{Mod}}}}
\rnc\mod[1]{\ensuremath{\mathop{\textup{mod-}#1}}\xspace}
\nc\Modtc{\Mod{\ct^c}}
\nc\pgldim[1]{\mathop{\rm pgldim}\,#1}
\nc\tstr{{\it t}--structure}
\nc\perf{^{\hbox{\rm\tiny perf}}}
\nc\op{^{\hbox{\rm\tiny op}}}
\nc\Th{^{\hbox{\rm\tiny th}}}
\nc\p{\varphi}
\nc\Dp{\{D\perf(\rs)\}}
\nc\ko{{K_0^{\hbox{\rm\tiny add}}}}

\begin{document}

\author{Andrew Ranicki}
\address{A.R. : School of Mathematics\newline
\indent University of Edinburgh \newline
\indent James Clerk Maxwell Building \newline
\indent King's Buildings \newline
\indent Mayfield Road \newline
\indent Edinburgh EH9 3JZ\newline
\indent SCOTLAND, UK}
\email{a.ranicki@ed.ac.uk}

\title[Noncommutative localization in algebraic $L$-theory]
{Noncommutative localization in algebraic $L$-theory}

\begin{abstract}\hskip2mm
Given a noncommutative (Cohn) localization $A \to \sigma^{-1}A$ which
is injective and stably flat we obtain a lifting theorem for
induced f.g. projective
$\rs$-module chain complexes and localization exact sequences in
algebraic $L$-theory, matching the algebraic $K$-theory localization
exact sequence of Neeman-Ranicki \cite{Neeman-Ranicki04} and
Neeman \cite{Neeman07}.
\end{abstract}


\maketitle

\tableofcontents

\section*{Introduction}
\label{S0}

The series of papers \cite{Neeman-Ranicki04}, \cite{Neeman07},
studied the algebraic $K$-theory of the
noncommutative (Cohn) localization $\sigma^{-1}A$ of a ring $A$ inverting
a collection $\sigma$ of morphisms of f.g.  projective left $A$-modules.
By definition, $\sigma^{-1}A$ is stably flat if
$${\rm Tor}^A_i(\sigma^{-1}A,\sigma^{-1}A)~=~0~~(i \geq 1)~.$$
An $(A,\sigma)$-module is an $A$-module $T$ which admits a f.g. projective
$A$-module resolution
$$0 \ri P \xymatrix{\ar[r]^s&} Q \ri T \ri 0$$
with $s:\sigma^{-1}P \to \sigma^{-1}Q$ an isomorphism of the
induced $\sigma^{-1}A$-modules.
For $A \ri \sigma^{-1}A$ which is injective and stably flat
we obtained an algebraic $K$-theory localization exact sequence
$$\dots \to K_n(A) \to K_n(\sigma^{-1}A) \to K_{n-1}(H(A,\sigma))
\to K_{n-1}(A) \to \dots$$
with $H(A,\sigma)$ the exact category of $(A,\sigma)$-modules.
\medskip

Let $C$ be a bounded $\sigma^{-1}A$-module chain complex such that each
$C_i=\sigma^{-1}P_i$ is induced from a f.g.  projective $A$-module
$P_i$.  The {\it chain complex lifting problem} is to decide if $C$ is
chain equivalent to $\sigma^{-1}D$ for a bounded chain complex $D$ of
f.g.  projective $A$-modules.  The problem has a trivial affirmative
solution for a commutative or Ore localization, by the clearing of
denominators, when $C$ is actually isomorphic to $\sigma^{-1}D$.  In
general, it is not possible to lift chain complexes: the injective
noncommutative localizations $A \to \sigma^{-1}A$ which are not stably
flat constructed in Neeman, Ranicki and Schofield \cite[Remark
2.13]{Neeman-Ranicki-Schofield} provide examples of induced f.g.
projective $\sigma^{-1}A$-module chain complexes of dimensions
$\geqslant 3$ which cannot be lifted.
\medskip

In \S1 we solve the chain complex lifting problem in the injective
stably flat case, obtaining the following results (Theorems
\ref{T1.3},\ref{T1.5})~:

\thm{0.1} For a stably flat injective noncommutative localization
$A \to \sigma^{-1}A$ every bounded chain
complex $C$ of induced f.g.  projective $\sigma^{-1}A$-modules is chain
equivalent to $\sigma^{-1}D$ for a bounded chain complex $D$ of f.g.
projective $A$-modules. Moreover, if $C$ is $n$-dimensional
$$C~:~\dots \to 0 \to C_n \to C_{n-1} \to \dots \to C_1 \to C_0 \to 0 \to \dots$$
then $D$ can be chosen to be $n$-dimensional.\hfill$\Box$
\ethm

In \S2 we consider the algebraic $L$-theory of a noncommutative
localization, obtaining the following results (Theorems
\ref{L2}, \ref{L3}, \ref{Lfin})~:

\thm{0.5} Let $A \ri  \sigma^{-1}A$ be a noncommutative localization
of a ring with involution $A$, such that $\sigma$ is invariant under
the involution.\\
{\rm (i)} There is a localization exact sequence of quadratic
$L$-groups
$$\xymatrix{\dots \ar[r] & L_n(A) \ar[r]&
L_n^I(\sigma^{-1}A) \ar[r]^-{\partial} &
L_n(A,\sigma)\ar[r] & L_{n-1}(A) \ar[r] & \dots}$$
with $I={\rm im}(K_0(A) \ri K_0(\sigma^{-1}A))$, and
$L_n(A,\sigma)$ the cobordism group of $\sigma^{-1}A$-contractible
$(n-1)$-dimensional quadratic Poincar\'e complexes over $A$.\\
{\rm (ii)} If $\sigma^{-1}A$ is stably flat over $A$
there is a localization exact sequence of symmetric $L$-groups
$$\xymatrix{\dots \ar[r] & L^n(A) \ar[r]&
L^n_I(\sigma^{-1}A) \ar[r]^-{\partial} &
L^n(A,\sigma)\ar[r] & L^{n-1}(A) \ar[r] & \dots}$$
with $L^n(A,\sigma)$ the cobordism group of
$\sigma^{-1}A$-contractible
$(n-1)$-dimensional symmetric Poincar\'e complexes over $A$.\\
{\rm (iii)} If $A \ri \sigma^{-1}A$ is injective then $L^n(A,\sigma)$
(resp. $L_n(A,\sigma)$) is the cobordism group of  $n$-dimensional
symmetric
(resp. quadratic) Poincar\'e complexes of $(A,\sigma)$-modules.
\hfill$\Box$
\ethm

The $L$-theory exact sequences of Theorem \ref{0.5} for an injective
Ore localization $A \ri \sigma^{-1}A$ (which is flat and hence stably
flat) were obtained in Ranicki \cite{Ranicki1981}.  The quadratic
$L$-theory exact sequence of \ref{0.5} (i) for arbitrary injective
$A \ri \sigma^{-1}A$ was obtained by Vogel \cite{Vogel1980}, \cite{Vogel1982}.
The symmetric $L$-theory exact sequence of \ref{0.5} (ii) is new.
\medskip

We refer to \cite{Ranicki1998,Ranicki2006} for some of the applications
of the algebraic $L$-theory of noncommutative localizations to topology.
\medskip

Amnon Neeman used to be a coauthor of the paper, but decided to
withdraw in May 2007.

\section{Lifting chain complexes}
\label{S1}

If $A\ri \rs$ is a stably flat localization, we know from
\cite[Theorem~0.4, Proposition~4.5 and Theorem~3.7]{Neeman-Ranicki04}
that the functor $Ti:\frac{D\perf(A)}{\car^c}\ri D\perf(\rs)$
is just an idempotent completion; it is fully faithful and all objects
in $D\perf(\rs)$ are, up to isomorphisms, direct summands of objects
in the image of $Ti$. A fairly easy consequence of this is the
following. Let $C\in D\perf(\rs)$ be the complex
\[
0 \ri \sm C^m\ri \sm C^{m+1}\ri\cdots\ri
\sm C^{n-1}\ri \sm C^n\ri0,
\]
with $C^i$ all finitely generated, projective $A$--modules.
Then there is complex $X\in D\perf(A)$
with $C\simeq \br\oti X$. That is, $C$ is homotopy equivalent to
the tensor product with $\rs$ of a perfect complex over the ring $A$.
In Section~\ref{S1} we prove this (Theorem~\ref{T1.3}), and then refine
the result to show that $X$ may be chosen to be a complex of the
form
\[
0 \ri X^m\ri X^{m+1}\ri\cdots\ri
X^{n-1}\ri X^n\ri0~.
\]
(Proof in Theorem~\ref{T1.5}).

\rmk{R0.1} The proof of Theorem~\ref{T1.3} relies
on the following fact about triangulated categories. Suppose $\ca$
is a full, triangulated subcategory of a triangulated category $\cb$,
and suppose all objects in $\cb$ are direct summands of objects of $\ca$.
An object $X\in\cb$ belongs to $\ca\subset\cb$ if and only if
$[X]\in K_0(\cb)$ lies in the image of $K_0(\ca)\ri K_0(\cb)$. This
fact may be found, for example, in~\cite[Proposition~4.5.11]{Neeman99},
but for the reader's convenience its proof is included here in Lemma~\ref{L1.1}
and Proposition~\ref{P1.2}.
\ermk

We begin by reminding the reader of some basic facts
about Grothendieck groups.
For any additive category $\ca$ we define
$\ko(\ca)$ to be the Grothendieck group of the split exact
category $\ca$. This means that the short exact sequences
in $\ca$ are precisely the split sequences.
It is well known that every element of $\ko(\ca)$
can be expressed as
\[
[X]-[Y]
\]
for $X$ and $Y$ objects of $\ca$. The expressions $[X]-[Y]$ and
$[X']-[Y']$ are equal in $\ko(\ca)$ if and only if there exists
an object $P\in\ca$ and an isomorphism
\[
X\oplus Y'\oplus P~=~ X'\oplus Y\oplus P.
\]

If $\ca$ happens to be a triangulated category, then $K_0(\ca)$
means the quotient of $\ko(\ca)$ by a subgroup
we will denote $T(\ca)$. The subgroup $T(\ca)$
is defined as the group generated by all
\[
[X]-[Y]+[Z],
\]
where there exists a distinguished triangle in $\ca$
\[
\CD
X @>>> Y @>>> Z @>>> \T X.
\endCD
\]
We prove:

\lem{L1.1}
Suppose $\cb$ is a triangulated category. Let $\ca$ be a full, triangulated
subcategory of $\cb$. Assume further that every object of $\cb$ is a direct
summand of an object in $\ca\subset\cb$.

Then the map $f:\ko(\ca)\ri \ko(\cb)$ induces a surjection
$T(\ca)\ri T(\cb)$. In symbols: $f\big(T(\ca)\big)=T(\cb)$.
\elem

\prf
Let $[X]-[Y]+[Z]$ be a generator of $T(\cb)\subset \ko(\cb)$.
We need to show it lies in the image of $T(\ca)\subset\ko(\ca)$.
Suppose therefore that
\[
\CD
X @>>> Y @>>> Z @>>> \T X
\endCD
\]
is a distinguished triangle in $\cb$. Because every object
of $\cb$ is a direct summand of an object in $\ca$, we can choose
objects $C$ and $D$ with
\[
X\oplus C,\qquad Z\oplus D
\]
both lying in $\ca$. But then we have a two
distinguished triangles in $\cb$
\[
\CD
X @>>> Y @>>> Z @>>> \T X \\
C @>>> C\oplus D @>>> D @>0>> \T C
\endCD
\]
and their direct sum is a distinguished triangle
\[
\CD
X\oplus C @>>> Y\oplus C\oplus D @>>> Z\oplus D @>>> \T (X\oplus C).
\endCD
\]
Two of the objects lie in $\ca$. Since the subcategory $\ca\subset\cb$
is full and triangulated, the entire distinguished
triangle lies in $\ca$. Thus
\[
[X\oplus C]-[Y\oplus C\oplus D]+[Z\oplus D]\quad=\quad
[X]-[Y]+[Z]
\]
lies in the image of $T(\ca)$.
\eprf

The next proposition is well-known; again, the proof is included for the
convenience of the reader.

\pro{P1.2}
Let the hypotheses be as in Lemma~\ref{L1.1}. That is,
suppose $\cb$ is a triangulated category. Let $\ca$ be a full, triangulated
subcategory of $\cb$. Assume further that every object of $\cb$ is a direct
summand of an object in $\ca\subset\cb$.

If $X$ is an object of $\cb$ and $[X]$ lies in the image of
the natural map $f:K_0(\ca)\ri K_0(\cb)$, then $X\in\ca$.
\epro

\prf
If we consider $[X]$ as an element of $\ko(\cb)$, then saying that
its image in $K_0(\cb)$ lies in the image of $K_0(\ca)\ri K_0(\cb)$
is equivalent to saying that, modulo $T(\cb)$, $[X]$ lies in the
image of $\ko(\ca)$. That is,
\[
[X]\in T(\cb)+ f\big(\ko(\ca)\big)\subset \ko(\cb).
\]
By Lemma~\ref{L1.1} we have that $f\big(T(\ca)\big)=T(\cb)$.
Thus
\[
\begin{array}{rcl}
 T(\cb)+f\big(\ko(\ca)\big)&=&f\big(T(\ca)\big)+ f\big(\ko(\ca)\big)\\
&=& f\big(\ko(\ca)\big).
\end{array}
\]
That means there exist objects $C$ and $D$ in $\ca\subset\cb$ and an identity
in $\ko(\cb)$
\[
[X]=[C]-[D].
\]
There must therefore be an object $P\in\cb$ and an isomorphism
\[
X\oplus D\oplus P~\simeq~C\oplus P.
\]
But $P$ is an object of $\cb$, hence a direct summand of an
object of $\ca$. There is an object $P'\in\cb$ with $P\oplus P'\in\ca$.
We have an isomorphism
\[
X\oplus D\oplus P\oplus P'~\simeq~C\oplus P\oplus P'.
\]
Putting $D'=D\oplus P\oplus P$ and $C'=C\oplus P\oplus P'$ we have
objects $C',D'$ in $\ca$, and a (split) distinguished triangle
\[
\CD
D' @>>> C' @>>> X @>>> \T D'.
\endCD
\]
Since $\ca\subset \cb$ is triangulated we conclude that $X\in\ca$.
\eprf

The relevance of these results to our work here is

\thm{T1.3}
Let $A\ri\rs$ be a stably flat localization of rings. Suppose
we are given a
perfect complex $C$ over $\rs$. Suppose further
that $C\in D\perf(\rs)$ is of the form
\[
0 \ri \sm C^m\ri \sm C^{m+1}\ri\cdots\ri
\sm C^{n-1}\ri \sm C^n\ri0
\]
where each $C^i$ is a finitely generated, projective $A$--module.
Then $C$ is homotopy equivalent to $\br\oti X$, for
some $X\in D\perf(A)$.
\ethm

\prf
The localization is stably flat. By~\cite[Theorem~0.4]{Neeman-Ranicki04}
the functor $T:\ct^c\ri D\perf(\rs)$ is an equivalence of categories.
By~\cite[Proposition~4.5 and Theorem~3.7]{Neeman-Ranicki04} we also
know that the functor $i:\frac{D\perf(A)}{\car^c}\ri \ct^c$
is fully faithful, and that every object in $\ct^c$ is
isomorphic to a direct summand
of an object in the image of $i$. Next we apply
Proposition~\ref{P1.2}, with $\cb=D\perf(\rs)$ and $\ca$ the
full subcategory
containing all objects isomorphic to $Ti(x)$, for any
$x\in \frac{D\perf(A)}{\car^c}$.

Now $C$ is an object of $D\perf(\rs)$, and in $K_0\big(
D\perf(\rs)
\big)$ we have an identity
\[
[C]=\sum_{\ell=-\infty}^\infty{(-1)}^\ell[\sm C^\ell]
\]
with
\[
[\sm C^\ell]=[\br\otimes_A^{}C^\ell]=[Ti C^\ell]
\]
certainly lying in the image of the map
\[
\CD
\ds
K_0(Ti):K_0\left(\frac
{D\perf(A)}{\car^c}
\right) @>>>
K_0\big(
D\perf(\rs)
\big).
\endCD
\]
Proposition~\ref{P1.2} therefore tells us that $C$ is
isomorphic to an object in the
image of the functor $Ti$. There exists a perfect complex
$X\in D\perf(A)$ and a homotopy equivalence $C\simeq\br\oti X$.
\eprf

The problem with Theorem~\ref{T1.3} is that it gives us no bound on
the length of the complex $X$ with $\br\oti X\simeq C$.
We really want to know

\thm{T1.5}
Let $A\ri\rs$ be a stably flat localization of rings.
Suppose $C\in D\perf(\rs)$
is the complex
\[
0 \ri \sm C^m\ri \sm C^{m+1}\ri\cdots\ri
\sm C^{n-1}\ri \sm C^n\ri0.
\]
Then the complex $X\in D\perf(A)$
with $C\simeq \br\oti X$, whose existence is guaranteed by
Theorem~\ref{T1.3}, may be chosen to be a complex
\[
0 \ri X^m\ri X^{m+1}\ri\cdots\ri
X^{n-1}\ri X^n\ri0~.
\]
\ethm

If $m=n$ this is easy. For $m<n$ we need to prove something. Our
proof will appeal to the results of~\cite[Section~4]{Neeman-Ranicki04}.
We remind the reader that this was the section which dealt with
the subcategories $\ck[m,n]$ of complexes in $\car^c$
vanishing outside the
range $[m,n]$. First we need a lemma.

\lem{Pbound.8}
Let $M$ and $N$ be any finitely generated projective $A$--modules.
We may view $M$ and $N$ as objects in the derived category $D\perf(A)$,
concentrated in degree $0$. Then any map in $\ct^c(\pi M,\pi N)$ can be
represented as ${\pi(\alpha)}^{-1}\pi(\beta)$, for some
$\alpha$, $\beta$ morphisms in $D\perf(A)$ as below
\[
\CD
M @>{\beta}>> Y @<{\alpha}<< N~.
\endCD
\]
The map $\alpha:N\ri Y$ fits in a triangle
\[
\CD
X @>>> N @>\alpha>> Y @>>> \Sigma X
\endCD
\]
and $X$ may be chosen to lie in $\ck[0,1]$.
\elem

\prf
By~\cite[Proposition~4.5 and Theorem~3.7]{Neeman-Ranicki04} we
know that the map
\[
\CD
i: \ds\frac
{D\perf(A)}{\car^c} @>>> \ct^c
\endCD
\]
is fully faithful. Therefore
\[
\ct^c(\pi M,\pi N)\eq
\frac
{D\perf(A)}{\car^c}(M,N).
\]
That is, any map $\pi M\ri\pi N$ can be written as
${\pi(\alpha)}^{-1}\pi(\beta)$, for some
$\alpha$, $\beta$ morphisms in $D\perf(A)$ as below
\[
\CD
M @>{\beta}>> Y @<{\alpha}<< N~.
\endCD
\]
The map $\alpha:N\ri Y$ fits in a triangle
\[
\CD
X @>>> N @>\alpha>> Y @>\beta>> \Sigma X
\endCD
\]
and $X$ may be chosen to lie in $\car^c$. What is not
clear is that we may choose $X$ in $\ck[0,1]\subset\car^c$.

The easy observation is that we may certainly modify
our choice of $X$ to lie
in $\ck\subset\car^c$. This follows from~%
\cite[Lemma~4.5]{Neeman07}, which tells us that
for any choice of $X$ as above there exists an $X'$ with $X\oplus X'$
isomorphic to an object in $\ck$. We have a distinguished
triangle
\[
\CD
X\oplus X' @>>> N
@>{\left(\begin{array}{c}\alpha \\0
\end{array}\right)
}>> Y\oplus \T X' @>\beta\oplus 1>> \Sigma (X\oplus X')
\endCD
\]
and a diagram
\[
\CD
M @>{\left(\begin{array}{c}\beta \\0
\end{array}\right)
}>> Y\oplus\T X' @<{\left(\begin{array}{c}\alpha \\0
\end{array}\right)
}<< N~,
\endCD
\]
and replacing our original choices by these we may assume
$X\in\ck$. Now we have to shorten $X$.

By~\cite[Lemma~4.7]{Neeman07}, there exists a triangle in
$\car^c$
\[
\CD
X' @>>> X @>>> X'' @>>> \Sigma X'
\endCD
\]
with $X'\in\ck[1,\infty)$ and $X''\in\ck(-\infty,1]$.
The composite $X'\ri X\ri N$ is a map from
$X' \in\ck[1,\infty)$ to $N\in\cs^{\leq0}$, which must
vanish. Hence we have that $X\ri N$ factors as
$X\ri X''\ri N$. We complete to a morphism of triangles
\[
\CD
X @>>> N @>\alpha>> Y @>>> \Sigma X \\
@VVV  @V1VV     @V\gamma VV         @VVV \\
X'' @>>> N @>\gamma\alpha>> Y'' @>>> \Sigma X''
\endCD
\]
and another representative of our morphism is the diagram
\[
\CD
M @>{\gamma\beta}>> Y'' @<{\gamma\alpha}<< N
\endCD
\]
We may, on replacing $Y$ by $Y''$, assume $X\in\ck(-\infty,1]$.

Applying~\cite[Lemma~4.7]{Neeman07}
again, we have that any
$X\in\ck(-\infty,1]$ admits a triangle
\[
\CD
X' @>>> X @>>> X'' @>>> \Sigma X'
\endCD
\]
with $X'\in\ck[0,1]$ and  $X''\in\ck(-\infty,0]$. Form the
octahedron
\[
\CD
X' @>>> N @>\alpha'>> Y' @>>> \Sigma X' \\
@VVV  @V1VV     @V\gamma VV         @VVV \\
X @>>> N @>\alpha>> Y @>>> \Sigma X \\
@. @. @VVV @VVV \\
@. @.   \Sigma X'' @>1>> \Sigma X''
\endCD
\]
The composite $M\ri Y\ri \Sigma X''$ is a map
from the projective module $M$, viewed as a complex
concentrated in degree 0, to $\Sigma X''\in\ck(\infty,-1]$.
This composite must vanish. The map
$\beta:M\ri Y$ therefore factors as
$M\stackrel{\beta'}\ri Y'
\stackrel\gamma\ri Y$, and our morphism in $\ct^c$ has a
representative
\[
\CD
M @>{\beta'}>> Y' @<{\alpha'}<< N
\endCD
\]
so that in the triangle
\[
\CD
X' @>>> N @>\alpha'>> Y' @>>> \Sigma X'
\endCD
\]
$X'$ may be chosen to lie in $\ck[0,1]$.
\eprf

Now we are ready for

\nin
{\bf Proof of Theorem~\ref{T1.5}.}\ \
We are given a complex $C\in D\perf(\rs)$ of the form
\[
0 \ri \sm C^m\ri \sm C^{m+1}\ri\cdots\ri
\sm C^{n-1}\ri \sm C^n\ri0.
\]
To eliminate the trivial case, assume $m\leq n+1$. Shifting,
we may assume $m=0$ and $n\geq1$.
Theorem~\ref{T1.3} guarantees that $C$ is homotopy equivalent to
$\br\oti D$, with $D\in D\perf(A)$. But $D$ need not be
supported on the interval $[0,n]$. We need to show how to
shorten $D$. Assume therefore that $D$ is supported on
$[-1,n]$. We will show how to replace $D$ by a complex supported
on $[0,n]$. Shortening a complex supported on $[0,n+1]$ is dual,
and we leave it to the reader.

We may
suppose therefore that
$D\in D\perf(A)$ is the complex
\[
\CD
\cdots @>>> 0 @>>> D^{-1} @>>> D^0 @>>> \cdots @>>> D^n@>>> 0 @>>> \cdots
\endCD
\]
and that there is a homotopy equivalence
of $\s^{-1}D$ with a shorter complex, that is a
commutative diagram
\[
\CD
 @>>> 0 @>>> \s^{-1}D^{-1} @>\partial>> \s^{-1}D^0
@>>> \cdots @>>> \s^{-1}D^n @>>> 0 @>>> \\
@. @VVV @VVV @VVV @. @VVV @VVV @. \\
 @>>> 0 @>>> 0 @>>> \sm C^0 @>>> \cdots @>>> \sm C^n @>>> 0 @>>> \\
@. @VVV @VVV @VVV @. @VVV @VVV @. \\
 @>>> 0 @>>> \s^{-1}D^{-1} @>\partial>> \s^{-1}D^0
@>>> \cdots @>>> \s^{-1}D^n @>>> 0 @>>>
\endCD
\]
so that the composite is homotopic to the identity.
In particular, there is a map $d:\s^{-1}D^0\ri \s^{-1}D^{-1}$
so that $d\partial:\s^{-1}D^{-1}\ri \s^{-1}D^{-1}$ is the
identity.

By~\cite[Proposition~3.1]{Neeman07} the map
$d:\s^{-1}D^0\ri \s^{-1}D^{-1}$
lifts uniquely to a map $d':\pi D^0\ri \pi D^{-1}$.
By Lemma~\ref{Pbound.8} the map $d'$ can be represented
as ${\pi(\alpha)}^{-1}\pi(\beta)$,
where $\alpha$ and $\beta$ are, respectively,
the chain maps
\[
\CD
 @>>> 0 @>>> 0 @>>> D^{-1} @>>> 0 @>>> \\
@. @VVV @VVV @VVV  @VVV @. \\
 @>>> 0 @>>> X @>r>> Y @>>> 0 @>>>
\endCD
\]
and
\[
\CD
 @>>> 0 @>>> 0 @>>> D^{0} @>>> 0 @>>> \\
@. @VVV @VVV @VgVV  @VVV @. \\
 @>>> 0 @>>> X @>r>> Y @>>> 0 @>>>
\endCD
\]
The fact that $\s^{-1}\alpha$ is an equivalence tells
us that the map $\s^{-1}r:\s^{-1}X\ri\s^{-1}Y$ is injective,
with cokernel $ \s^{-1}D^{-1}$. The fact that $\alpha^{-1}\beta$
agrees with $d'$ means that the composite
\[
\CD
\s^{-1}D^{0} @>\s^{-1}g>> \s^{-1}Y @>>> \text{\rm Coker}(\s^{-1}r)
\endCD
\]
is just the map $d:\s^{-1}D^{0}\ri \s^{-1}D^{-1}$.
Let $X$ be the chain complex
\[
\CD
@>>> 0  @>>> D^0 \oplus X@>{\begin{pmatrix}
\partial & 0 \\ g & r
\end{pmatrix}}>> D^{1}\oplus Y @>>>
\cdots @>>> D^n @>>> 0 @>>>
\endCD
\]
Let $f:X\ri D$
be the
natural map of chain complexes
\[
\CD
\ri 0@>>> 0  @>>> D^0 \oplus X@>{\begin{pmatrix}
\partial & 0 \\ g & r
\end{pmatrix}}>> D^{1}\oplus Y @>>>
\cdots @>>> D^n \ri 0 \ri \\
@. @VVV @V{\pi_1^{}}VV @VV{\pi_1^{}}V @. @VVV @. \\
\ri 0 @>>> D^{-1} @>>> D^0 @>>\partial> D^{1} @>>>
\cdots @>>> D^n \ri 0 \ri
\endCD
\]
where the vertical maps labelled $\pi_1^{}$ are the projections to
the first factor of the direct sum. The map $\s^{-1}f$ is easily
seen to be homotopy equivalence. Thus $\s^{-1}X$ is homotopy
equivalent to $\s^{-1}D\cong C$.\hfill{$\Box$}

\section{Algebraic $L$-theory}
\label{Ltheory}

An {\it involution} on a ring $A$ is an anti-automorphism
$$A \ri A~;~r \mapsto \overline{r}~.$$
The involution is used to regard a left $A$-module $M$ as a right
$A$-module by
$$M \times A \ri M ~;~(x,r) \mapsto \overline{r} x~.$$
The {\it dual} of a (left) $A$-module $M$ is the $A$-module
$$M^*~=~{\rm Hom}_A(M,A)~,~A \times  M^* \ri M^*~;~
(r,f) \mapsto (x \mapsto f(x)\overline{r})~.$$
The {\it dual} of an $A$-module morphism $s:P \ri Q$ is the $A$-module
morphism
$$s^*~:~Q^* \ri P^*~;~f \mapsto (x \mapsto f(s(x)))~.$$
If $M$ is f.g. projective then so is $M^*$, and
$$M \ri M^{**}~;~x \mapsto (f \mapsto \overline{f(x)})$$
is an isomorphism which is used to identify $M^{**}=M$.

\begin{hypothesis} \label{linj1}
In this section, we assume that
\begin{itemize}
\item[(i)] $A$ is a ring with involution,
\item[(ii)] the duals of morphisms $s:P \ri Q$ in $\sigma$
are morphisms $s^*:Q^* \ri P^*$ in $\sigma$,
\item[(iii)] $\epsilon \in A$ is a central unit such that
$\overline{\epsilon}=\epsilon^{-1}$ {\rm (}e.g. $\epsilon= \pm 1${\rm )}.
\end{itemize}
The noncommutative localization $\sigma^{-1}A$
is then also a ring with involution, with $\epsilon \in \sigma^{-1}A$
a central unit such that $\overline{\epsilon}=\epsilon^{-1}$.
\hfill$\Box$
\end{hypothesis}

We review briefly the chain complex construction of the f.g. projective
$\epsilon$-quadratic $L$-groups $L_*(A,\epsilon)$
and the $\epsilon$-symmetric $L$-groups $L^*(A,\epsilon)$.
Given an $A$-module chain complex $C$ let
the generator $T \in \bZ_2$ act on the $\bZ$-module
chain complex $C\otimes_AC$ by the $\epsilon$-transposition duality
$$T_{\epsilon}~:~C_p \otimes_A C_q \ri C_q \otimes_A C_p~:~x \otimes y
\mapsto (-1)^{pq}\epsilon y \otimes x~.$$
Let $W$ be the standard free $\bZ[\bZ_2]$-module resolution of
$\mathbb Z$
$$W~:~\dots \ri \bZ[\bZ_2]
\xrightarrow[]{1-T} \bZ[\bZ_2] \xrightarrow[]{1+T} \bZ[\bZ_2]
\xrightarrow[]{1-T} \bZ[\bZ_2]~.$$
The {\it $\epsilon$-symmetric} (resp. {\it $\epsilon$-quadratic})
{\it $Q$-groups} of $C$ are the
$\bZ_2$-hypercohomology (resp. $\bZ_2$-hyperhomology)
groups of $C \otimes_AC$
$$\begin{array}{l}
Q^n(C,\epsilon)~=~H^n(\bZ_2;C\otimes_AC)~=~
H_n({\rm Hom}_{\bZ[\bZ_2]}(W,C\otimes_AC))~,\\[1ex]
Q_n(C,\epsilon)~=~H_n(\bZ_2;C\otimes_AC)~=~
H_n(W\otimes_{\bZ[\bZ_2]}(C\otimes_AC))~.
\end{array}$$
The $Q$-groups are chain homotopy invariants of $C$.
There are defined forgetful maps
$$\begin{array}{l}
1+T_{\epsilon}~:~Q_n(C,\epsilon) \ri Q^n(C,\epsilon)~;~
\psi \mapsto (1+T_{\epsilon})\psi~,\\[1ex]
Q^n(C,\epsilon) \ri H_n(C \otimes_A C)~;~\phi
\mapsto \phi_0~.
\end{array}$$
For f.g. projective $C$ the function
$$C\otimes_AC \ri {\rm Hom}_A(C^*,C)~;~
x \otimes y \mapsto (f \mapsto \overline{f(x)}y)$$
is an isomorphism of $\bZ[\bZ_2]$-module chain complexes, with
$T \in \bZ_2$ acting on ${\rm Hom}_A(C^*,C)$ by $\theta \mapsto \epsilon \theta^*$.
The element $\phi_0 \in H_n(C\otimes_AC)=H_n({\rm Hom}_A(C^*,C))$ is
a chain homotopy class of $A$-module chain maps $\phi_0:C^{n-*} \ri C$.
\medskip

An {\it $n$-dimensional $\epsilon$-symmetric complex over $A$} $(C,\phi)$
is a bounded f.g. projective $A$-module chain complex $C$ together with
an element $\phi \in Q^n(C,\epsilon)$. The complex $(C,\phi)$
is {\it Poincar\'e} if the $A$-module chain map $\phi_0:C^{n-*} \ri C$
is a chain equivalence.

\exm{form}
A 0-dimensional $\epsilon$-symmetric Poincar\'e complex $(C,\phi)$ over $A$
is essentially the same as a nonsingular $\epsilon$-symmetric form
$(M,\lambda)$ over $(A,\sigma)$, with $M=(C_0)^*$ a f.g. projective
$A$-module and
$$\lambda~=~\phi_0~:~M \times M \ri A$$
a sesquilinear pairing such that the adjoint
$$M \ri M^*~;~x \mapsto (y \mapsto \lambda(x,y))$$
is an $A$-module isomorphism.
\eexm

See pp.\ 210--211 of \cite{Ranicki1998} for the notion of an
{\it $\epsilon$-symmetric (Poincar\'e) pair}.
The {\it boundary} of an $n$-dimensional
$\epsilon$-symmetric complex $(C,\phi)$ is the $(n-1)$-dimensional
$\epsilon$-symmetric Poincar\'e complex
$$\partial (C,\phi)~=~(\partial C,\partial \phi)$$
with $\partial C=C(\phi_0:C^{n-*} \ri C)_{*+1}$ and $\partial \phi$
as defined on p.\ 218 of \cite{Ranicki1998}.
The {\it $n$-dimensional $\epsilon$-symmetric $L$-group}
$L^n(A,\epsilon)$ is the cobordism group of $n$-dimensional
$\epsilon$-symmetric Poincar\'e complexes $(C,\phi)$ over $A$
with $C$ $n$-dimensional.
In particular, $L^0(A,\epsilon)$ is the Witt group of nonsingular
$\epsilon$-symmetric forms over $A$.
\medskip

An $n$-dimensional $\epsilon$-symmetric complex $(C,\phi)$ over $A$
is {\it $\sigma^{-1}A$-Poincar\'e} if the $\sigma^{-1}A$-module chain map
$\sigma^{-1}\phi_0:\sigma^{-1}C^{n-*} \ri \sigma^{-1}C$
is a chain equivalence, in which case $\sigma^{-1}(C,\phi)$ is an
$n$-dimensional $\epsilon$-symmetric Poincar\'e complex over
$\sigma^{-1}A$.
\medskip

The {\it $n$-dimensional $\epsilon$-symmetric $\Gamma$-group}
$\Gamma^n(A\ri \sigma^{-1}A,\epsilon)$ is the cobordism group of
$n$-dimensional $\epsilon$-symmetric $\sigma^{-1}A$-Poincar\'e
complexes $(C,\phi)$ over $A$ such that $\sigma^{-1}C$ is chain
equivalent to an $n$-dimensional induced f.g.  projective
$\sigma^{-1}A$-module chain complex.  The {\it $n$-dimensional
$\epsilon$-symmetric $L$-group} $L^n(A,\sigma,\epsilon)$ is the
cobordism group of $(n-1)$-dimensional $\epsilon$-symmetric Poincar\'e
complexes over $A$ $(C,\phi)$ such that $C$ is
$\sigma^{-1}A$-contractible, i.e.  $\sigma^{-1}C \simeq 0$.  \medskip

Similarly in the $\epsilon$-quadratic case, with groups
$L_n(A,\epsilon)$, $\Gamma_n(A \ri \sigma^{-1}A,\epsilon)$,
$L_n(A,\sigma,\epsilon)$. The $\epsilon$-quadratic $L$- and $\Gamma$-groups
are 4-periodic
$$\begin{array}{l}
L_n(A,\epsilon)~=~L_{n+2}(A,-\epsilon)~=~L_{n+4}(A,\epsilon)~,\\[1ex]
\Gamma_n(A\ri \sigma^{-1}A,\epsilon)~=~
\Gamma_{n+2}(A\ri \sigma^{-1}A,-\epsilon)~=~\Gamma_{n+4}(A\ri \sigma^{-1}A,\epsilon)~,\\[1ex]
L_n(A,\sigma,\epsilon)~=~L_{n+2}(A,\sigma,-\epsilon)~=~L_{n+4}(A,\sigma,\epsilon)~.
\end{array}$$

\pro{L1} For any ring with involution $A$ and noncommutative localization
$\sigma^{-1}A$ there is defined a localization exact sequence
of $\epsilon$-symmetric $L$-groups
$$\xymatrix{
\dots \ar[r] & L^n(A,\epsilon) \ar[r]&
\Gamma^n(A \ri \sigma^{-1}A,\epsilon) \ar[r]^-{\partial} &
L^n(A,\sigma,\epsilon)\ar[r] & L^{n-1}(A,\epsilon) \ar[r] & \dots~.}$$
Similarly in the $\epsilon$-quadratic case, with an exact sequence
$$\xymatrix{\dots \ar[r] & L_n(A,\epsilon) \ar[r]&
\Gamma_n(A \ri \sigma^{-1}A,\epsilon) \ar[r]^-{\partial} &
L_n(A,\sigma,\epsilon)\ar[r] & L_{n-1}(A,\epsilon) \ar[r] & \dots~.}$$
\epro
\prf The relative group of $L^n(A,\epsilon) \ri \Gamma^n(A \ri
\sigma^{-1}A,\epsilon)$ is the cobordism group of
$n$-dimensional $\epsilon$-symmetric $\sigma^{-1}A$-Poincar\'e
pairs over $A$ $(f:C \ri D,(\delta\phi,\phi))$ with $(C,\phi)$ Poincar\'e.
The effect of algebraic surgery on $(C,\phi)$ using this pair is a
cobordant $(n-1)$-dimensional $\epsilon$-symmetric Poincar\'e complex
$(C',\phi')$ with $C'$ $\sigma^{-1}A$-contractible. The function
$(f:C \ri D,(\delta\phi,\phi)) \mapsto (C',\phi')$ defines an isomorphism
between the relative group and $L^n(A,\sigma,\epsilon)$.
\eprf

Define
$$I~=~{\rm im}(K_0(A)\ri K_0(\sigma^{-1}A))~,$$
the subgroup of $K_0(\sigma^{-1}A)$ consisting of the projective classes
of the f.g.  projective $\sigma^{-1}A$-modules induced from f.g.  projective
$A$-modules. By definition, $L^n_I(\sigma^{-1}A,\epsilon)$ is the cobordism
group of $n$-dimensional $\epsilon$-symmetric Poincar\'e complexes over
$\sigma^{-1}A$
$(B,\theta)$ such that $[B] \in I$.
There are evident morphisms of $\Gamma$- and $L$-groups
$$\begin{array}{l}
\sigma^{-1}\Gamma^*~:~\Gamma^n(A \ri \sigma^{-1}A,\epsilon) \ri
L_I^n(\sigma^{-1}A,\epsilon)~;~(C,\phi) \mapsto \sigma^{-1}(C,\phi)~,\\[1ex]
\sigma^{-1}\Gamma_*~:~\Gamma_n(A \ri \sigma^{-1}A,\epsilon) \ri
L^I_n(\sigma^{-1}A,\epsilon)~;~(C,\psi) \mapsto \sigma^{-1}(C,\psi)~.
\end{array}$$
In general, the morphisms $\sigma^{-1}\Gamma^*,\sigma^{-1}\Gamma_*$
need not be isomorphisms, since a bounded f.g.  projective
$\sigma^{-1}A$-module chain complex $D$ with $[D] \in I$ need not be
chain equivalent to $\sigma^{-1}C$ for a bounded f.g.  projective
$A$-module chain complex $C$.
\medskip

It was proved in Chapter 3 of Ranicki \cite{Ranicki1981} that
if $A \ri \sigma^{-1}A$ is an injective Ore localization then the morphisms
$\sigma^{-1}Q^*,\sigma^{-1}Q_*$,
$\sigma^{-1}\Gamma^*,\sigma^{-1}\Gamma_*$ are isomorphisms,
so that there are defined localization exact sequences
for both the $\epsilon$-symmetric and the $\epsilon$-quadratic $L$-groups
$$\xymatrix@R-23pt{
\dots \ar[r] & L^n(A,\epsilon) \ar[r]&
L^n_I(\sigma^{-1}A,\epsilon) \ar[r]^-{\partial} &
L^n(A,\sigma,\epsilon)\ar[r] & L^{n-1}(A,\epsilon) \ar[r] & \dots~,\\
\dots \ar[r] & L_n(A,\epsilon) \ar[r]&
L^I_n(\sigma^{-1}A,\epsilon) \ar[r]^-{\partial} &
L_n(A,\sigma,\epsilon)\ar[r] & L_{n-1}(A,\epsilon) \ar[r] & \dots~.}$$
Special cases of these sequences were obtained by Milnor-Husemoller,
Karoubi, Pardon, Smith, Carlsson-Milgram.
\medskip

Let $G\pi:D(A) \to D(A)$ be the functor of Proposition 6.1 of \cite{Neeman-Ranicki04},
with $D(A)$ the derived category of $A$.
For any bounded f.g. projective $A$-module chain complex $C$
the natural $A$-module chain map
$$\mathop{\varinjlim}\limits_{(B,\beta)} B ~=~G\pi(C) \ri \sigma^{-1}C$$
induces morphisms
$$\begin{array}{l}
\sigma^{-1}Q^*~:~\mathop{\varinjlim}\limits_{(B,\beta)}
Q^n(B,\epsilon)~=~Q^n(G\pi(C),\epsilon) \ri Q^n(\sigma^{-1}C,\epsilon)~,\\[1ex]
\sigma^{-1}Q_*~:~\mathop{\varinjlim}\limits_{(B,\beta)}
Q_n(B,\epsilon)~=~Q_n(G\pi(C),\epsilon) \ri Q_n(\sigma^{-1}C,\epsilon)
\end{array} $$
with the direct limits taken over all the bounded f.g.  projective
$A$-module chain complexes $B$ with a chain map $\beta:C\ri B$ such
that $\sigma^{-1}\beta:\sigma^{-1}C \ri \sigma^{-1}B$ is a
$\sigma^{-1}A$-module chain equivalence.  The natural projection
$D\otimes_AD \ri D\otimes_{\sigma^{-1}A}D$ is an isomorphism for any
bounded f.g.  projective $\sigma^{-1}A$-module chain complex $D$ (since
this is already the case for $D=\sigma^{-1}A$), so the $Q$-groups of
$\sigma^{-1}C$ are the same whether $\sigma^{-1}C$ is regarded as an
$A$-module or $\sigma^{-1}A$-module chain complex.

\thm{L2}  {\rm (Vogel \cite{Vogel1982}, Theorem 8.4)}
For any ring with involution $A$ and noncommutative localization
$\sigma^{-1}A$ the morphisms
$$\sigma^{-1}\Gamma_*~:~\Gamma_n(A \ri \sigma^{-1}A,\epsilon) \ri
L^I_n(\sigma^{-1}A,\epsilon)~;~(C,\psi) \mapsto \sigma^{-1}(C,\psi)$$
are isomorphisms, and there is a localization exact sequence
of $\epsilon$-quadratic $L$-groups
$$\xymatrix{
\dots \ar[r] & L_n(A,\epsilon) \ar[r]&
L^I_n(\sigma^{-1}A,\epsilon) \ar[r]^-{\partial} &
L_n(A,\sigma,\epsilon)\ar[r] & L_{n-1}(A,\epsilon) \ar[r] & \dots~.}$$
\ethm
\prf By algebraic surgery below the middle dimension it suffices
to consider only the special cases $n=0,1$. In effect, it was proved
in \cite{Vogel1982} that $\sigma^{-1}Q_*$ is an isomorphism for
0- and 1-dimensional $C$.
\eprf

It was claimed in Proposition 25.4 of Ranicki \cite{Ranicki1998}
that $\sigma^{-1}\Gamma^*$ is also an isomorphism, assuming (incorrectly)
that the chain complex lifting problem can always be solved.
However, we do have :

\thm{L3}
If $\sigma^{-1}A$ is a noncommutative localization of a ring with
involution $A$ which is stably flat over $A$, there is a localization
exact sequence of $\epsilon$-symmetric $L$-groups
$$\xymatrix{
\dots \ar[r] & L^n(A,\epsilon) \ar[r]&
L^n_I(\sigma^{-1}A,\epsilon) \ar[r]^-{\partial} &
L^n(A,\sigma,\epsilon)\ar[r] & L^{n-1}(A,\epsilon) \ar[r] & \dots~.}$$
\ethm
\prf For any bounded f.g.  projective $A$-module chain complex $C$ the
natural $A$-module chain map $G\pi(C) \ri \sigma^{-1}C$ induces
isomorphisms in homology
$$H_*(G\pi(C))~\cong~H_*(\sigma^{-1}C)~.$$
Thus the natural $\bZ[\bZ_2]$-module chain map
$$G\pi(C)\otimes_AG\pi(C) \ri \sigma^{-1}C\otimes_A\sigma^{-1}C~=~
\sigma^{-1}C\otimes_{\sigma^{-1}A}\sigma^{-1}C$$
induces isomorphisms of $\epsilon$-symmetric $Q$-groups
$$\sigma^{-1}Q^*~:~\mathop{\varinjlim}\limits_{(B,\beta)}
Q^n(B,\epsilon) \ri Q^n(\sigma^{-1}C,\epsilon)$$
(and also isomorphisms $\sigma^{-1}Q_*$ of $\epsilon$-quadratic $Q$-groups).
By Theorem \ref{0.1} every $n$-dimensional induced f.g. projective
$\sigma^{-1}A$-module chain complex $D$ is chain equivalent to $\sigma^{-1}C$
for an $n$-dimensional f.g. projective $A$-module chain complex $C$, with
$$Q^n(D,\epsilon)~=~Q^n(\sigma^{-1}C,\epsilon)~=~
\mathop{\varinjlim}\limits_{(B,\beta)}Q^n(B,\epsilon) ~.$$
It follows that the morphisms of $\epsilon$-symmetric $\Gamma$- and $L$-groups
$$\sigma^{-1}\Gamma^*~:~\Gamma^n(A \ri \sigma^{-1}A,\epsilon) \ri
L^n_I(\sigma^{-1}A,\epsilon)~;~(C,\phi) \mapsto \sigma^{-1}(C,\phi)$$
are also isomorphisms, and the localization exact sequence is given
by Proposition \ref{L1}.
\eprf

\begin{hypothesis} \label{linj2}
For the remainder of this section, we assume Hypothesis \ref{linj1}
and also that $A \ri \sigma^{-1}A$ is an injection.\hfill$\Box$
\end{hypothesis}
As in Proposition 2.2 of \cite{Neeman07}
it follows that all the morphisms in $\s$ are injections.
\medskip

We shall now  generalize the results of Ranicki
\cite{Ranicki1981} and Vogel \cite{Vogel1980} to prove that under
Hypotheses \ref{linj1},\ref{linj2} the relative $L$-groups
$L^*(A,\sigma,\epsilon)$, $L_*(A,\sigma,\epsilon)$ in the $L$-theory
localization exact sequences are the $L$-groups of $H(A,\sigma)$ with
respect to the following duality involution.
\medskip

Define the {\it torsion dual}
of an $(A,\sigma)$-module $M$ to be the $(A,\s)$-module
$$M\widehat{~}~=~{\rm Ext}^1_A(M,A)~,$$
using the involution on $A$ to define the left $A$-module
structure. If $M$ has f.g. projective $A$-module resolution
$$0 \ri P_1 \xrightarrow[]{s} P_0 \ri M \ri 0$$
with $s \in \sigma$ the torsion dual $M\widehat{~}$ has the dual f.g.
projective $A$-module resolution
$$0 \ri P_0^* \xrightarrow[]{s^*} P_1^* \ri M\widehat{~} \ri 0$$
with $s^* \in \sigma$.
\medskip

\pro{iso} Let $M={\rm coker}(s:P_1 \ri P_0)$,
$N={\rm coker}(t:Q_1 \ri Q_0)$ be $(A,\sigma)$-modules.\\
{\rm (i)} The adjoint of the pairing
$$M \times M\widehat{~} \ri \sigma^{-1}A/A~;~(g \in P_0,f\in P_1^*)
\mapsto fs^{-1}g$$
defines a natural $A$-module isomorphism
$$M\widehat{~} \ri {\rm Hom}_A(M,\sigma^{-1}A/A)~;~f \mapsto (g \mapsto fs^{-1}g)~.$$
{\rm (ii)} The natural $A$-module morphism
$$M \ri M\widehat{~}\widehat{~}~;~x \mapsto (f \mapsto \overline{f(x)})$$
is an isomorphism.\\
{\rm (iii)} There are natural identifications
$$\begin{array}{l}
M \otimes_A N~=~{\rm Tor}^A_0(M,N)~=~{\rm Ext}^1_A(M\widehat{~},N)~=~
H_0(P\otimes_AQ)~,\\[1ex]
{\rm Hom}_A(M\widehat{~},N)~=~{\rm Tor}^A_1(M,N)~=~
{\rm Ext}^0_A(M\widehat{~},N)~=~H_1(P\otimes_AQ)~.
\end{array}$$
The functions
$$\begin{array}{l}
M \otimes_A N \ri N \otimes_A M~;~x \otimes y \mapsto y \otimes x~,\\[1ex]
{\rm Hom}_A(M\widehat{~},N) \ri {\rm Hom}_A(N\widehat{~},M)~;~
f \mapsto f\widehat{~}
\end{array}$$
determine transposition isomorphisms
$$T~:~{\rm Tor}_i^A(M,N) \ri {\rm Tor}_i^A(N,M)~~(i=0,1)~.$$
{\rm (iv)} For any finite subset
$V=\{v_1,v_2,\dots,v_k\} \subset M\otimes_AN$ there exists an exact sequence of
$(A,\sigma)$-modules
$$0 \ri N \ri L  \ri \oplus_k M\widehat{~} \ri 0$$
such that $V \subset {\rm ker}(M\otimes_AN \ri M\otimes_AL)$.
\epro
\prf (i)  Apply the snake lemma to the morphism of short exact sequences
$$\xymatrix{
0 \ar[r] &{\rm Hom}_A(P_0,A) \ar[r] \ar[d]^-{\displaystyle{s^*_{}}} &
{\rm Hom}_A(P_0,\sigma^{-1}A) \ar[r] \ar[d]^-{\displaystyle{s^*_1}} &
{\rm Hom}_A(P_0,\sigma^{-1}A/A) \ar[r] \ar[d]^-{\displaystyle{s^*_2}}& 0 \\
0 \ar[r] &{\rm Hom}_A(P_1,A) \ar[r] &
{\rm Hom}_A(P_1,\sigma^{-1}A) \ar[r] &
{\rm Hom}_A(P_1,\sigma^{-1}A/A) \ar[r] & 0}$$
with $s^*$ injective, $s^*_1$ an isomorphism and $s^*_2$ surjective,
to verify that the $A$-module morphism
$$M\widehat{~}~=~{\rm coker}(s^*) \ri
{\rm Hom}_A(M,\sigma^{-1}A/A)~=~{\rm ker}(s^*_2)$$
is an isomorphism.\\
(ii) Immediate from the identification
$$s^{**}~=~s~:~(P_0)^{**}~=~P_0 \ri (P_1)^{**}~=~P_1~.$$
(iii) Exercise for the reader.\\
(iv) Lift each $v_i \in M\otimes_AN$ to an element
$$v_i \in P_0 \otimes_AQ_0~=~{\rm Hom}_A(P_0^*,Q_0)~~(1 \leq i \leq k)~.$$
The $A$-module morphism defined by
$$u~=~\begin{pmatrix}s^* & 0 & 0 & \dots & 0 \\
0 & s^* & 0 & \dots & 0 \\
0 & 0 & s^* & \dots & 0 \\
\vdots & \vdots & \vdots &\ddots & \vdots \\
v_1 & v_2 & v_3 & \dots & t \end{pmatrix}~:~
U_1~=~(\oplus_k P_0^*) \oplus
Q_1 \ri U_0~=~(\oplus_k P_1^*)\oplus Q_0$$
is in $\sigma$, so that $L={\rm coker}(u)$ is an $(A,\sigma)$-module with
a f.g. projective $A$-module resolution
$$0 \ri U_1 \xrightarrow[]{u} U_0 \ri L \ri 0~.$$
The short exact sequence of 1-dimensional f.g. projective $A$-module chain
complexes
$$0 \ri Q \ri U \ri \oplus_k P^{1-*} \ri 0$$
is a resolution of a short exact sequence of $(A,\sigma)$-modules
$$0 \ri N \ri L  \ri \oplus_k M\widehat{~} \ri 0~.$$
The first morphism in the exact sequence
$${\rm Tor}^A_1(M,\oplus_k M\widehat{~}) \ri M\otimes_AN \ri
M\otimes_AL \ri M \otimes_A(\oplus_k M\widehat{~}) \ri 0$$
sends $1_i \in {\rm Tor}^A_1(M,\oplus_kM\widehat{~})=
\oplus_k {\rm Hom}_A(M\widehat{~},M\widehat{~})$
to $v_i \in {\rm ker}(M\otimes_AN \ri M\otimes_AL)$.
\eprf

Given an  $(A,\sigma)$-module chain complex $C$ define the {\it
$\epsilon$-symmetric} (resp.  {\it $\epsilon$-quadratic}) {\it torsion
$Q$-groups} of $C$ to be the $\bZ_2$-hypercohomology (resp.
$\bZ_2$-hyperhomology) groups of the $\epsilon$-transposition
involution $T_{\epsilon}=\epsilon T$ on the $\bZ$-module chain complex
${\rm Tor}_1^A(C,C)={\rm Hom}_A(C\widehat{~},C)$
$$\begin{array}{l}
Q_{\rm tor}^n(C,\epsilon)~=~H^n(\bZ_2;{\rm Tor}^A_1(C,C))~=~
H_n({\rm Hom}_{\bZ[\bZ_2]}(W,{\rm Tor}^A_1(C,C)))~,\\[1ex]
Q^{\rm tor}_n(C,\epsilon)~=~H_n(\bZ_2;{\rm Tor}^A_1(C,C))~=~
H_n(W\otimes_{\bZ[\bZ_2]}({\rm Tor}^A_1(C,C)))~.
\end{array}$$
There are defined forgetful maps
$$\begin{array}{l}
1+T_{\epsilon}~:~Q^{\rm tor}_n(C,\epsilon) \ri Q^n_{\rm tor}(C,\epsilon)~;~
\psi \mapsto (1+T_{\epsilon})\psi~,\\[1ex]
Q_{\rm tor}^n(C,\epsilon) \ri H_n({\rm Tor}^A_1(C,C))~;~\phi
\mapsto \phi_0~.
\end{array}$$
The element $\phi_0 \in H_n({\rm Tor}^A_1(C,C))$ is
a chain homotopy class of $A$-module chain maps
$\phi_0:C^{n-}\widehat{~} \ri C$.
\medskip

An {\it $n$-dimensional $\epsilon$-symmetric complex over $(A,\sigma)$}
$(C,\phi)$ is a bounded $(A,\sigma)$-module chain complex $C$ together with
an element $\phi \in Q_{\rm tor}^n(C,\epsilon)$. The complex $(C,\phi)$
is {\it Poincar\'e} if the $A$-module chain maps
$\phi_0:C^{n-}\widehat{} \ri C$ are chain equivalences.

\exm{linking form}
A 0-dimensional $\epsilon$-symmetric Poincar\'e complex $(C,\phi)$ over $(A,\sigma)$
is essentially the same as a nonsingular $\epsilon$-symmetric linking form
$(M,\lambda)$ over $(A,\sigma)$, with $M=(C_0)\widehat{~}$ an
$(A,\s)$-module and
$$\lambda~=~\phi_0~:~M \times M \ri \sigma^{-1}A/A$$
a sesquilinear pairing such that the adjoint
$$M \ri M\widehat{~}~;~x \mapsto (y \mapsto \lambda(x,y))$$
is an $A$-module isomorphism.
\eexm

The {\it $n$-dimensional torsion $\epsilon$-symmetric $L$-group}
$L_{\rm tor}^n(A,\sigma,\epsilon)$ is the cobordism group of
$n$-dimensional $\epsilon$-symmetric Poincar\'e complexes $(C,\phi)$ over
$(A,\sigma)$, with $C$ $n$-dimensional.  In particular,
$L^0_{\rm tor}(A,\sigma,\epsilon)$ is the Witt group of nonsingular
$\epsilon$-symmetric linking forms over $(A,\sigma)$.
\medskip

Similarly in the $\epsilon$-quadratic case, with torsion $L$-groups
$L^{\rm tor}_n(A,\sigma,\epsilon)$.  The $\epsilon$-quadratic torsion
$L$-groups are 4-periodic
$$L^{\rm tor}_n(A,\sigma,\epsilon)~=~
L^{\rm tor}_{n+2}(A,\sigma,-\epsilon)~=~L^{\rm tor}_{n+4}(A,\sigma,\epsilon)~.$$

\thm{Lfin}
If $A \ri \sigma^{-1}A$ is injective the relative $L$-groups in the localization exact
sequences of Proposition \ref{L1}
$$\xymatrix@R-23pt{
\dots \ar[r] & L^n(A,\epsilon) \ar[r]&
\Gamma^n(A \ri \sigma^{-1}A,\epsilon) \ar[r]^-{\partial} &
L^n(A,\sigma,\epsilon)\ar[r] & L^{n-1}(A,\epsilon) \ar[r] & \dots\\
\dots \ar[r] & L_n(A,\epsilon) \ar[r]&
\Gamma_n(A \ri \sigma^{-1}A,\epsilon) \ar[r]^-{\partial} &
L_n(A,\sigma,\epsilon)\ar[r] & L_{n-1}(A,\epsilon) \ar[r] & \dots}
$$
are the torsion $L$-groups
$$\begin{array}{l}
L^*(A,\sigma,\epsilon)~=~L_{\rm tor}^*(A,\sigma,\epsilon)~,\\[1ex]
L_*(A,\sigma,\epsilon)~=~L^{\rm tor}_*(A,\sigma,\epsilon)~.
\end{array}$$
\ethm
\prf For any bounded $(A,\sigma)$-module chain complex $T$
there exists a bounded f.g. projective $A$-module chain complex $C$
with a homology equivalence $C \ri T$. Working as in \cite{Vogel1980}
there is defined a distinguished triangle of $\bZ[\bZ_2]$-module chain complexes
$$\Sigma {\rm Tor}^A_1(T,T) \ri
C\otimes_AC \ri T\otimes_AT \ri \Sigma^2{\rm Tor}^A_1(T,T)$$
with $\bZ_2$ acting by the $\epsilon$-transposition $T_{\epsilon}$
on the $\bZ$-module chain complex ${\rm Tor}^A_1(T,T)$
and by the $(-\epsilon)$-transpositions $T_{-\epsilon}$
on $C\otimes_AC$ and $T\otimes_AT$, inducing long exact sequences
$$\xymatrix@R-23pt{
\dots \ar[r] & Q^n_{\rm tor}(T,\epsilon) \ar[r]&
Q^{n+1}(C,-\epsilon) \ar[r] &
Q^{n+1}(T,-\epsilon) \ar[r] & Q^{n-1}_{\rm tor}(T,\epsilon) \ar[r] & \dots\\
\dots \ar[r] & Q_n^{\rm tor}(T,\epsilon) \ar[r]&
Q_{n+1}(C,-\epsilon) \ar[r] &
Q_{n+1}(T,-\epsilon) \ar[r] & Q_{n-1}^{\rm tor}(T,\epsilon) \ar[r] & \dots~.}
$$
Passing to the direct limits over all the bounded $(A,\sigma)$-module
chain complexes $U$ with a homology equivalence $\beta:T \ri U$ use
Proposition \ref{iso} (iv) to obtain
$$\begin{array}{l}
\mathop{\varinjlim}\limits_{(U,\beta)} Q^{n+1}(U,-\epsilon)~=~0~,\\[1ex]
\mathop{\varinjlim}\limits_{(U,\beta)} Q_{n+1}(U,-\epsilon)~=~0
\end{array}$$
and hence
$$\begin{array}{l}
\mathop{\varinjlim}\limits_{(U,\beta)}
Q^n_{\rm tor}(U,\epsilon)~=~Q^{n+1}(C,-\epsilon)~,\\[1ex]
\mathop{\varinjlim}\limits_{(U,\beta)}
Q_n^{\rm tor}(U,\epsilon)~=~Q_{n+1}(C,-\epsilon)~.
\end{array}$$
\eprf

\rmk{final} The identification $L_*(A,\sigma,\epsilon)=L^{\rm
tor}_*(A,\sigma,\epsilon)$ for noncommutative $\sigma^{-1}A$ was first
obtained by Vogel \cite{Vogel1980}.
\ermk

\ifx\undefined\bysame
\newcommand{\bysame}{\leavevmode\hbox to3em{\hrulefill}\,}
\fi

\end{document}